\documentclass[11pt]{amsart}
\usepackage[hypertex]{hyperref}
\usepackage{mathrsfs}
\usepackage{amssymb}
\usepackage{amsfonts}
\usepackage{graphicx}
\usepackage{subfigure}
\usepackage{mathptmx}
\usepackage{latexsym,amsmath,amssymb,amsfonts,amsthm}
\usepackage{amsbsy}
\newcommand{\LA}[1]{\mbox{\LARGE $#1$}}

\newtheorem{theorem}{Theorem}[section]

\newtheorem{proposition}[theorem]{Proposition}

\newtheorem{definition}[theorem]{Definition}
\title{ The Isoperimetric Problem in Higher Codimension}
\author{Frank Morgan  \and  Isabel M.C.\ Salavessa}
\date{\today}
\thanks{This work has been partially supported by the National Science
Foundation and the Funda\c{c}\~{a}o para a Ci\^{e}ncia e a Tecnologia,
programms  PTDC/MAT/101007/2008, and  PTDC/MAT/118682/2010.}
\address{
  Department of Mathematics and Statistics,
 Bronfman Science Center,
Williams College, 18 Hoxsey Street,
Williamstown, MA 01267, USA}\email{Frank.Morgan@williams.edu}
\address{ Centro de F\'{\i}sica das Interac\c{c}\~{o}es
Fundamentais, Instituto Superior T\'{e}cnico,  Technical University
of Lisbon, Edif\'{\i}cio Ci\^{e}ncia, Piso 3, Av.\ Rovisco Pais,
1049-001 Lisboa, Portugal}\email{isabel.salavessa@ist.utl.pt}
\begin{document}
\maketitle
\noindent
\begin{abstract}
We consider three generalizations of the isoperimetric problem to higher
codimension and provide results on equilibrium, stability, and minimization.
\end{abstract}
\section{Introduction}
\noindent
The classical isoperimetric problem in an $n$-dimensional Riemannian manifold
seeks an $(n-1)$-dimensional surface $S$ of least area bounding a region $R$ 
of prescribed  volume. To generalize the problem to $m$-dimensional surfaces
$S$ ($1 \leq  m \leq  n-2$) requires a notion of enclosed volume. We present
three alternatives:\\[3mm]
$(1)$ infimum $\mathrm{v}(S)$ of volumes of $(m+1)$-dimensional surfaces 
bounded by $S$,\\[3mm]
$(2)$ $\omega$-volume $\int_S\omega $  for some given smooth $m$-form $\omega$,
 \\[3mm]
$(3)$ in $\mathbf{R}^n$ multi-volume, i.e., volume enclosed by projection to 
each axis $(m+1)$-dimensional vector subspace of $\mathbf{R}^n$, or 
equivalently prescribed $\omega$-volume for all $m$-forms $\omega$ with 
$\mathrm{d}\omega$ constant.\\

 For the first notion, perhaps the most natural,  Almgren (\cite{Alm1}, 1986)
 proved that in $\mathbf{R}^n$, round spheres are uniquely isoperimetric.

        The second notion was introduced by Salavessa (\cite{S1}, 2010), 
actually in terms of $\Omega = \mathrm{d}\omega$; note that for any surface 
$R$ bounded  by $S$,
$$\int_R \Omega  = \int_S \omega.$$
\noindent
Given an exact form $\Omega$, $\omega$ is defined up to a closed form. 
If $\Omega$ is a constant $(m+1)$-form in $\mathbf{R}^n$, it follows from 
Almgren's  result that round spheres are uniquely isoperimetric. 
Salavessa \cite{S1} proves 
that round spheres uniquely satisfy some strong stability hypotheses.

The third notion was introduced by Morgan (\cite{M4}, 2000), who characterized
isoperimetric curves (not necessarily round) and gave examples of non-round
isoperimetric surfaces.

\begin{figure}[hbt]
\begin{tabular}{cc}
%\resizebox{!}{220pt}{\includegraphics{Fig1a.pdf}}
\resizebox{!}{220pt}{\includegraphics{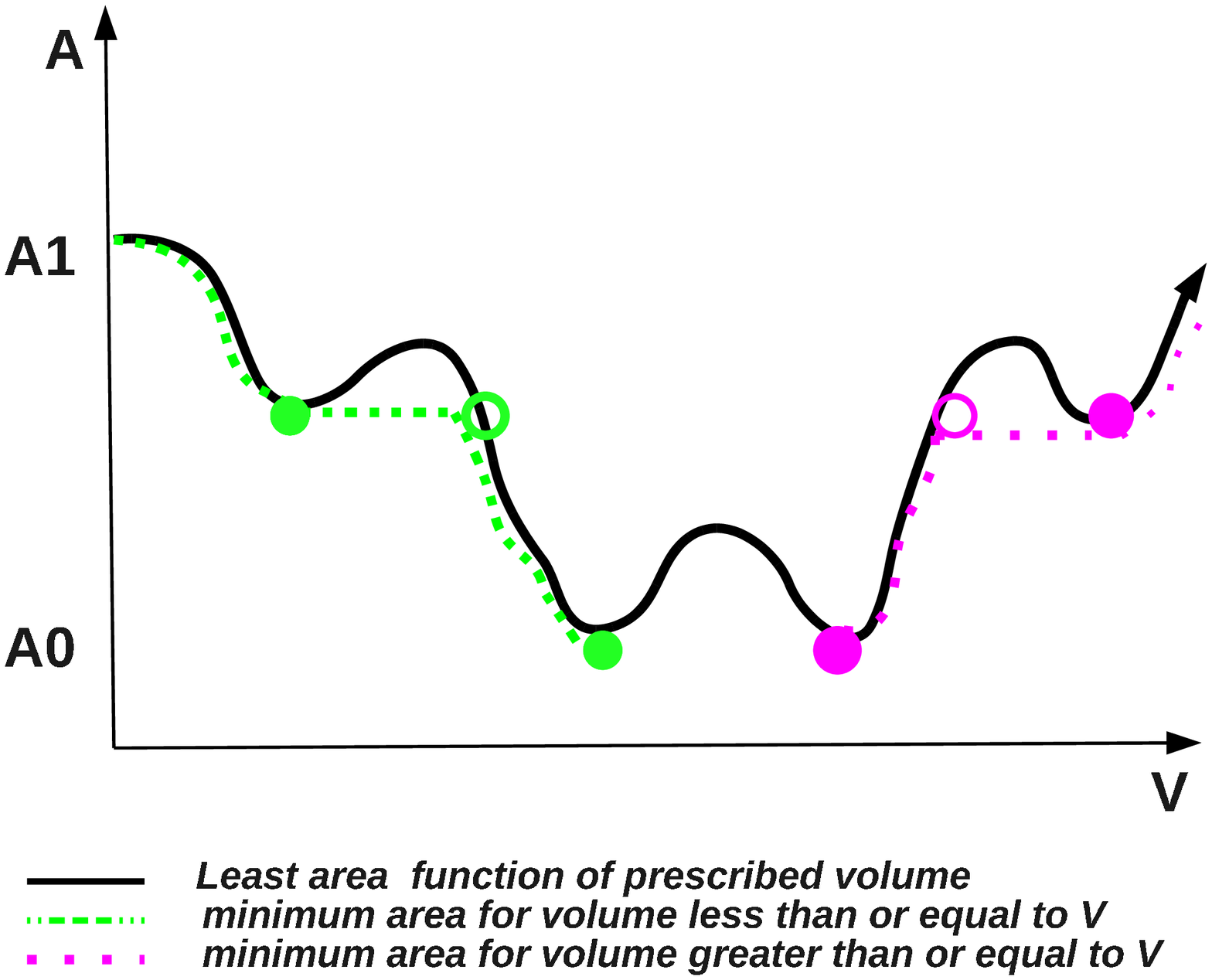}}
&\quad
%\resizebox{!}{220pt}{\includegraphics{Fig1b.pdf}}
\resizebox{!}{220pt}{\includegraphics{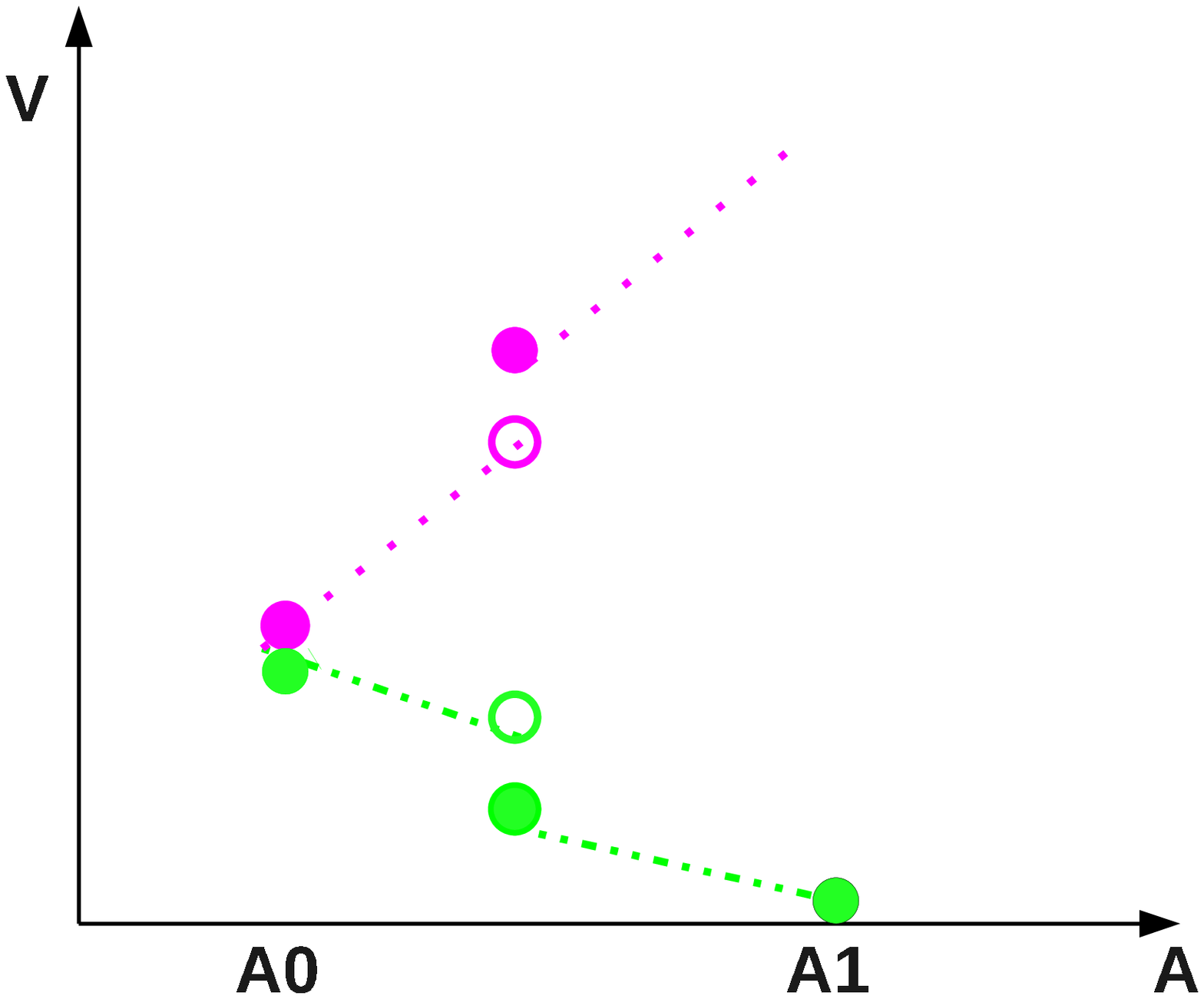}}
\end{tabular}\\[-13mm]
\caption{\footnotesize 
Least area (in black) is a continuous 
function of  prescribed  volume, but minimum and  maximum
 volume are not continuous  functions of  prescribed area.
In green  is  minimum area for volume less than or equal to $V$. 
In purple is minimum area for volume greater than or equal to $V$.
}
\label{fig1ab}
\end{figure}

We could more generally consider surfaces $S$ with prescribed boundary as well
as prescribed volume. In case (1), volume must then be measured with respect
to a given reference surface with the given boundary other than itself. 
In the closely related (higher dimensional) thread problem 
(see \cite{E}, [\cite{D}, Chap.\ 10], \cite{N,M2}),
which makes sense only for given reference surface, one fixes an area 
smaller than the reference surface
and minimizes volume. The fact that one is minimizing rather than maximizing 
volume makes it easy to prove that the thread has constant mean curvature 
inside the volume-minimizing surface where it is smooth; otherwise one could 
move the surface inward where the curvature is small and outward less where 
it is larger, preserving area but reducing volume, because volume is at most
the volume of the perturbed surface. A further difficulty for our case of 
prescribed volume is that there is no obvious perturbation preserving volume. 
The thread problem is roughly equivalent to minimizing area for a prescribed 
upper bound on volume. Least area is a continuous function of prescribed 
volume, but unless it is decreasing-increasing, minimum and maximum volume 
are not continuous functions of prescribed area. 
Figure~\ref{fig1ab} suggests possible relationships between area and volume, 
although we do not know a specific example that exhibits all these 
possibilities.

In case (1), we could also work in the larger context of unoriented surfaces.

To allow our surfaces to have singularities, we work in the context of the 
locally integral currents of geometric measure theory [\cite{M1}, 
Chaps.\ 4 and 9].

In earlier work R.\ Gulliver \cite{G1,G2}, 
F.\ Duzaar and M.\ Fuchs  (\cite{DF1, DF3},  and especially 
[\cite{DF2}, Thm.\ 3.2]), and Duzaar and 
K.\ Steffen \cite{DS}, seek surfaces 
with prescribed mean curvature vector by minimizing $A-\lambda V$. 
Gulliver [\cite{G1}, p.\ 118] gives one interpretation of a helical minimizer 
as  the  path of ``a charged particle moving in a magnetic field.''

This paper provides a unified treatment on minimizing area for the three 
notions (1)-(3) of prescribed volume. Section 2 discusses equilibrium 
conditions. Section 3 discusses existence and regularity. Section 4 considers 
the question of whether round spheres are the only isoperimetric or stable 
surfaces. We conjecture that in $\mathbf{R}^n$, round $m$-spheres 
$S_0$ are the only 
smooth stable surfaces $S$ for given volume $\mathrm{v}(S)$ or given $\Omega$-volume 
for constant $\Omega$ (although not for given multi-volume).

\section{Stationary Surfaces}
\noindent
This section presents the equilibrium conditions for the isoperimetric problem 
for the three types of volume constraints, generalizing the equilibrium 
condition of constant mean curvature of codimension 1. Higher codimension 
presents new issues of smoothness and degeneracy.

We will consider perturbations $S_t$ of an $m$-dimensional surface $S = S_0$
under nice smooth families $F_t$ $(0 \leq  t < t_1)$ of diffeomorphisms of 
$M$, with $F_0$ the identity. More specifically, we will assume that $F_t$
 is $C^3$ with spatial derivatives of orders 2 and 3 bounded, which includes 
scaling in $\mathbf{R}^n$. If $S$ has infinite area, we will assume that the 
$F_t$ equal the identity outside a fixed compact set. 

\begin{definition}
We call $S$  stationary  if the (one-sided) first derivative of the 
area $A = |S|$ of $S$ is nonnegative  whenever $S_t$ respect the volume 
constraint. We call $S$  stable  if small perturbations respecting 
the volume constraint have no less area. 
\end{definition}
\noindent
If $S$ is stationary, then the second variation depends only on the initial
variation vectorfield $\mathbf{v} = \partial F\!/\! \partial t. $

In the classical case of codimension 1 $(m = n-1)$, a stationary surface has
\em generalized mean curvature \em  $\mathbf{H}$ (defined almost everywhere) 
of constant  magnitude and normal to the surface, i.e., for every smooth 
variation  vectorfield $\mathbf{v}$, initially

$$ \mathrm{d}A\!/\!\mathrm{d}\,t = - \int_S (n-1) \mathbf{H} \cdot \mathbf{v}.$$
\noindent
This case is easy because volume $V$ also varies smoothly and non-degenerately; 
initially
$$ \mathrm{d}V\!\!/\!\mathrm{d}\,t = - \int_S \mathbf{n}\cdot \mathbf{v} ,$$
where $\mathbf{n}$ is the inward unit normal, defined almost everywhere. 
In higher codimension $V$ need not vary smoothly, as when $S$  bounds multiple 
volume-minimizing surfaces. Nor need $V$ vary non-degenerately: smooth families
 $F_t$ with $\mathrm{d}V\!\!/\!\mathrm{d}\,t$  initially $0$ sometimes cannot
be modified to keep $V$ constant (see \S 4.2 and \S 4.3). 
The generalized mean curvature
vector exists as long as $\mathrm{d}A\!/\!\mathrm{d}\,t$ is a bounded operator; 
see Allard \cite{A1} 
for details in a general (``varifold'') setting.

The rest of this section attempts to recover a constant-magnitude generalized
mean curvature vector for stationary surfaces in higher codimension for the 
three definitions of prescribed volume. In the most difficult first case of
prescribed volume $\mathrm{v}(S)$, Proposition 2.2 requires a strong smoothness 
hypothesis, while the more useful Proposition 2.3 uses a stronger notion of 
\em stationary. \em  In the other two cases of prescribed $\omega$-volume and 
prescribed multi-volume, volume varies smoothly but  degeneracy can be an issue.

\begin{proposition} Let $S$ be a boundary in a smooth 
Riemannian manifold $M$. Suppose that\\[2mm]
$(1)$ $~~$ \parbox[t]{130mm}{there is a nonzero measurable vectorfield 
$\mathbf{G}$ on $S$ such that  for any smooth variation vectorfield 
$\mathbf{v}$, the volume $V = \mathrm{v}(S)$  is smooth and initially 
$\mathrm{d}V\!\!/\!\mathrm{d}\,t = 
-\int_S \mathbf{G}\cdot \mathbf{v}$.}\\[2mm]
Then $S$ is stationary if and only if the generalized mean curvature vector 
$\mathbf{H}$ is a constant times $\mathbf{G}$.
\end{proposition}
\noindent 
{\em Remarks.}  If $S$ bounds a volume-minimizing surface $R$, then any such 
$\mathbf{G}$ is weakly the inward unit conormal, i.e., 
$\mathrm{d}|R|/\mathrm{d}\,t$ is initially 
$-\int_S \mathbf{G}\cdot \mathbf{v}$. If two such volume-minimizing
surfaces or sufficiently regular minimal surfaces are indecomposable and smooth
submanifolds with boundary at some point of $S$, they are equal, as follows by 
partial differential equations [\cite{M3}, 
Sect.\ 7], using for continuation the 
indecomposability of $S$ and the fact that volume-minimizing hypersurfaces are 
regular except possibly for a codimension-2 singular set \cite{Alm2}. 
In particular, 
if $S$ is a smooth, connected submanifold of $\mathbf{R}^n$,  then the 
volume-minimizing surface is unique because volume-minimizing surfaces are 
regular at extreme points of $S$ by Allard's boundary regularity 
theorem \cite{A2}. Conversely,\\[3mm]
 (2)$~~$ \parbox[t]{130mm}{ \em  we conjecture that hypothesis $(1)$ holds 
whenever  $S$ bounds a unique volume-minimizing surface. \em  }\\[3mm]
On the other hand, (1) clearly  fails if $S$ bounds two different 
volume-minimizing surfaces with distinct  conormals. For an extreme 
negative example, let $M$ be the round 2-sphere and let $S$ be two 
antipodal points.

By work of B.\ White \cite{Wh}, (2) holds in $\mathbf{R}^n$ and in compact 
real-analytic $n$-dimensional ambients $M$, as long as all volume-minimizing 
surfaces are smoothly immersed manifolds with boundary, which in turn holds 
if $S$ has dimension $n- 2 \leq  5$ or if $S$ has dimension $1$ and we admit 
unoriented volume-minimizing surfaces  ([\cite{M1}, Chap.\ 8], \cite{A2}; 
unoriented 2-dimensional area-minimizing surfaces have 
no branch points or singularities except where sheets cross orthogonally). 
Furthermore almost every $S$ bounds a unique volume-minimizing surface 
([\cite{M3}, Thm.\ 7.1 and Rmk.] with \cite{Alm2}).

Apparently many smooth surfaces with (even parallel) mean curvature vector of
constant length are not stationary for prescribed volume, such as smooth 
minimal submanifolds $S$ of the unit sphere in $\mathbf{R}^n$ 
(at least for $n \leq  7$)  bounding unique volume-minimizing surfaces other 
than the cone,  such as $\mathbf{S}^1(1/\sqrt{5}) 
\times \mathbf{S}^2(2/\sqrt{5})$ in 
$\mathbf{R}^5$,  if indeed that bounds a unique volume-minimizing surface. 
(If $S$ were stationary, the volume-minimizing surface, smooth along $S$
 by Allard's boundary regularity theorem, would by Proposition 2.2 have 
radially inward conormal and therefore equal the cone by the PDE
argument described earlier in these remarks.) 
Conversely, we doubt that all stationary surfaces have parallel mean curvature 
and give a probable counterexample in a manifold in the Remarks after 
Proposition 2.3. If one allows prescribed boundary as well as
prescribed volume, Yau's characterization of 2-dimensional surfaces with 
parallel mean curvature (Sect.\ 4.4) implies that some 2-dimensional 
isoperimetric surfaces in $\mathbf{R}^4$ have nonparallel mean curvature, 
namely when the boundary is not contained in some $\mathbf{S}^3$ or 
$\mathbf{R}^3$ and the surface is non-minimal (as it must be for large 
prescribed volume).

One example where (1) holds is a round $m$-sphere $S$ in $\mathbf{R}^n$, 
with $\mathbf{G}$ the inward unit conormal to the flat ball and $\mathbf{H}$
 proportional to $\mathbf{G}$; here an easy lower bound on volume is
provided by the projection into the $(m+1)$-plane containing $S$. 
Conjecture (2) would imply that one nonround example is 
$S = \mathbf{S}^3\times \mathbf{S}^3$ in $\mathbf{R}^8$, 
with $\mathbf{G}$ the inward unit conormal to
the cone over $S$, which is famously volume-minimizing, 
but we don't see how to obtain the requisite lower bound on volume. 
Proposition 2.4 provides an  alternative proof that $S$ is stationary.\\[4mm]
{\em Proof of Proposition $2.2$.} 
 If $S$ bounds a volume-minimizing surface $R$, since 
$\mathrm{d}V \leq  \mathrm{d}|R|$ for
$\mathrm{d}|R|$ positive or negative, by smoothness 
$\mathrm{d}V = \mathrm{d}|R|$, 
$\mathrm{d}|R|/\mathrm{d}\,t = 
\mathrm{d}V\!\!/\!\mathrm{d}\,t = -\int_S\mathbf{G}\cdot \mathbf{v}$, 
and $|\mathbf{G}| \leq  1$.

Suppose that the generalized mean curvature vector $\mathbf{H}$ 
exists and is a constant times $\mathbf{G}$. Then if $V$ is constant, 
initially $\mathrm{d}V\!\!/\!\mathrm{d}\,t = -\int_S \mathbf{G}\cdot 
\mathbf{v} = 0$, 
so $\mathrm{d}A\!/\!\mathrm{d}\,t 
= -\int_S \mathbf{H}\cdot \mathbf{v} = 0$, i.e., $S$ is stationary.

Conversely, suppose that $S$ stationary. Since $\mathbf{G}$ is nonzero, we can 
choose variation vectorfields $\mathbf{v}_1$ and $\mathbf{v}_2$ with 
disjoint supports such that initially $\mathrm{d}V\!\!/\!\mathrm{d}\,t$ 
is nonzero
for each of them. Given a point $p$ of $S$, consider a neighborhood of 
$p$ disjoint from the support of $\mathbf{v}_1$ or $\mathbf{v}_2$, 
say $\mathbf{v}_1$, and let $\mathbf{v}$ be a smooth variation vectorfield 
supported in that neighborhood. Now some linear combination 
$\mathbf{w} = \mathbf{v} + c_1\mathbf{v}_1$  has 
$\mathrm{d}V\!\!/\!\mathrm{d}\,t = 0$
 and hence by smoothness comes from a two-sided $(-t_1 < t < t_1)$
 volume-preserving family of diffeomorphisms of the form
$$F_t(x) = \mathrm{exp}_x(  t\mathbf{v} + \varphi(t) c_1\mathbf{v}_1),$$
with $\varphi(0)=0$ and $\varphi' (0)=1$.
Since $S$ is stationary, initially $\mathrm{d}A\!/\!\mathrm{d}\,t = 0$. 
Consequently if
$\mathrm{d}V\!\!/\!\mathrm{d}\,t$ for $\mathbf{v}$ is nonzero, 
then $\mathrm{d}A\!/\!\mathrm{d}V$ 
for $\mathbf{v}$
 is the same as it is for $\mathbf{v}_1$. In particular, 
$\mathrm{d}A\!/\!\mathrm{d}V$ for
$\mathbf{v}_2$ is the same as it is for $\mathbf{v}_1$. 
Therefore $\mathrm{d}A\!/\!\mathrm{d}V$ is a constant $c$. 
It follows that $c\,\mathbf{G}$ is a
generalized mean curvature vector for $S$.

\begin{proposition}  Let $S$ be a boundary in $\mathbf{R}^n$ with finite area 
and volume $\mathrm{v}(S)$. Let $H_0 =|S|/(m+1)\mathrm{v}(S)$. 
If $\mathrm{d}A\!/\!\mathrm{d}\,t 
\geq  0$ under smooth 
families of diffeomorphisms $F_t$ for which $\Delta \mathrm{v}(S)\geq  0$, 
then $S$ has generalized mean curvature $\mathbf{H}$ of magnitude bounded 
by $H_0$. 
Conversely if $S$ is smooth with mean curvature vector $\mathbf{H}$
 a constant nonnegative multiple of the inward
unit conormal of a volume-minimizing surface, then 
$\mathrm{d}A\!/\!\mathrm{d}\,t 
\geq  0$
 under smooth families of
diffeomorphisms $F_t$ for which $\Delta \mathrm{v}(S) \geq 0$, and 
$|\mathbf{H}| = H_0$.
\end{proposition} 
\noindent
\em Proof.  \em  The technical difficulty is that although 
the area of $S = S_0$  varies smoothly under smooth perturbations $S_t$, 
the volume $\mathrm{v}(S_t)$ may not. Nevertheless for any smooth variation 
vectorfield  $\mathbf{v}$, the change in $V\!(t) = \mathrm{v}(S_t)$ satisfies
\begin{equation}
 \Delta V \geq - |\Delta t| \int_S |\mathbf{v}| - o(\Delta t)
\end{equation}
because if there were surfaces bounded by $S_t$ of smaller volume, 
adding on the volume swept out by the $S_t$ would yield a surface bounded by 
$S_0$ of less volume than $V\!(0)$, a contradiction. 
Given $\epsilon > 0$,  consider rescalings by a factor $1 + a\,t$
 with $a\,$ chosen such that initially
$$ \mathrm{d}V\!\!/\!\mathrm{d}\,t = 
a\,(m+1)\,V(0) = \int_S |\mathbf{v}|  +  \epsilon$$
and hence
$$ \mathrm{d}A\!/\!\mathrm{d}\,t = a\, m \,A(0) = m \,H_0 \,
\left(\int_S |\mathbf{v}|  +  \epsilon \right).$$
After combining the 
original family with such rescalings,
(1) becomes 
$\Delta V \geq \epsilon \Delta t - o(\Delta t) \geq 0$ for 
$0 \leq  t < t_1$ and hence by hypothesis 
$\mathrm{d}A\!/\!\mathrm{d}\,t \geq  0$. 
Therefore for the original family 
$$\mathrm{d}A\!/\!\mathrm{d}\,t \geq -m\,H_0\left(\int_S|\mathbf{v}| 
+ \epsilon\right);$$ 
since this holds for all $\epsilon > 0$, 
$~\mathrm{d}A\!/\!\mathrm{d}t \geq - m\, H_0 \int_S |\mathbf{v}|$. 
Since this holds also 
 for $-\mathbf{v}$, $|\mathrm{d}A\!/\!\mathrm{d}\,t| \leq  m\,H_0$, 
and this holds for every 
smooth variation $\mathbf{v}$. It follows that there is a generalized 
curvature vector $\mathbf{H}$ of magnitude bounded by  $H_0$.

Conversely, suppose that $S$ is smooth with mean curvature vector 
$\mathbf{H}$ a constant nonnegative multiple of the inward unit conormal 
$\mathbf{n}$ of a  volume-minimizing surface $R$.
For any smooth family of diffeomorphisms for which 
$\Delta \mathrm{v}(S) \geq  0$,
$$ 0 \leq  \mathrm{d}|R|/\mathrm{d}\,t = -\int_S \mathbf{n}\cdot \mathbf{v}.$$
Since $\mathbf{H}$  is a constant positive multiple of $\mathbf{n}$, 
$\mathrm{d}A\!/\!\mathrm{d}\, t = -\int_S m\, 
\mathbf{H}\cdot \mathbf{v} \geq 0$. 

Under scaling, $\mathrm{d}A\!/\!\mathrm{d}\,t
= -\int_S  m \,\mathbf{H}\cdot \mathbf{v}$  and 
$\mathrm{d}V\!\!/\!\mathrm{d}\,t = -\int_S \mathbf{n}\cdot \mathbf{v}$. 
Also for scaling  $A^{1/m}=c V^{1/(m+1)}$ and
 $\mathrm{d}A\!/\!\mathrm{d}V =  m\, A\!/\!(m+1)\, V = m\, H_0$. 
Therefore $\,\mathbf{H} = H_0\mathbf{n}$.
\\[4mm]
{\em Remarks. \em }
 The first statement of the converse and its proof hold in any smooth
Riemannian manifold. The hypothesis, when $R$ is not smooth along $S$, 
need only hold weakly: $\mathrm{d}|R|/\mathrm{d}\,t$ proportional to $-\int_S
\mathbf{H}\cdot \mathbf{v}.$

In $\mathbf{R}^2 \times [0,\epsilon]$ with the top and bottom identified with 
an appropriate slight twist, the helix is stationary (by the converse, 
assuming that the helicoid is area minimizing) and probably isoperimetric. 

\begin{proposition}
For $p \geq  1$, $\mathbf{S}^p \times \mathbf{S}^p$ in $\mathbf{R}^{2p+2}$
 is strongly stationary as in Proposition $2.3$.
\end{proposition}
\noindent
\em Proof.  \em  For $p \geq 3$ the cone is famously volume minimizing 
(\cite{BDG}, see [\cite{M1}, \S 10.7]), and the result follows 
immediately from  Proposition 2.3. For the general case let $R_1$ be a 
volume-minimizing surface bounded by $S$. By Allard's boundary regularity 
theorem \cite{A2}, $R_1$ is a smooth submanifold with boundary along $S$. 
Let $R_2$ be its image  under the symmetry switching the first two 
coordinates with the last two  coordinates. Then $\mathbf{H}$, 
which is in the unique symmetric normal direction to $
\mathbf{S}^p \times \mathbf{S}^p$, must be 
proportional to the sum $\mathbf{n}_1 + \mathbf{n}_2$
 of the conormals. For any smooth family of
diffeomorphisms for which $\Delta \mathrm{v}(S) \geq  0$,
$$ 0 \leq  \mathrm{d}|R_i|/\mathrm{d}\,t = -\int_S \mathbf{n}_i \cdot \mathbf{v}.$$ 
Since $\mathbf{H}$ is a proportional to $\mathbf{n}_1 + \mathbf{n}_2$, 
$\,\mathrm{d}A\!/\!\mathrm{d}\,t = -\int_S  m\, \mathbf{H}\cdot \mathbf{v} 
\geq  0$.\\[1mm]

\noindent
We now consider the second case of prescribed $\omega$-volume. Salavessa
 [\cite{S1}, Thm.\ 2.1] proves the following equilibrium condition in 
the narrower 
context where $\mathrm{d}\omega  \LA{\llcorner} T$  has
constant length.
\begin{proposition}[cf.\ \cite{S1},  Thm.\ 2.1] Consider a cycle 
(surface without boundary) S with unit tangent $m$-vector $T$ 
in a smooth Riemannian manifold $M$ with smooth $m$-form $\omega$. $S$ is
stationary for positive prescribed $\omega$-volume if the mean curvature 
vector $\mathbf{H}$ is proportional to 
$\mathrm{d}\omega\LA{ \llcorner} T$ (weakly).  Further 
suppose that $\mathrm{d}\omega \LA{\llcorner} T$ is not identically 0 or that
$M = \mathbf{R}^n$ and $\mathrm{d}\omega$ is constant. 
If $S$ is stationary, then 
the mean curvature vector $\mathbf{H}$ is
proportional to $\mathrm{d}\omega \LA{\llcorner} T$ (weakly).

If $M = \mathbf{R}^n$, $\mathrm{d}\omega$  is constant and simple, and $S$
 is stationary, smooth, connected, and bounded, then $S$ lies in an 
associated $(m+1)$-plane and is round.
\end{proposition}
\noindent
That the mean curvature vector $\mathbf{H}$ is proportional to 
$\mathrm{d}\omega  \LA{\llcorner} T$ (weakly) means that there is a
constant $c$ such that for any smooth variation vectorfield $\mathbf{v}$, 
initially
$$ \mathrm{d}A\!/\!\mathrm{d}\,t 
= c \int_S (\mathrm{d}\omega \LA{\llcorner} T) (\mathbf{v}) = c \int_S 
\mathrm{d}\omega (T \!\wedge \mathbf{v}).$$
(The generalized mean curvature $\mathbf{H}$  is characterized by 
$\mathrm{d}A\!/\!\mathrm{d}\,t 
= - \int_S m\, \mathbf{H}\cdot \mathbf{v}$ for any smooth
variation vectorfield $\mathbf{v}$, and one often identifies the vector 
$\mathbf{H}$ with the 1-form $\mathbf{H}\cdot $.)\\[3mm]
\noindent
\em Remark. \em The additional hypothesis for the converse is necessary. 
For example, let $S$ be any embedding of the hypersphere of finite area in 
$\mathbf{R}^n$, with inside $U$ and outside $V$. Let
$f, g$ be nonnegative $C^{\infty}$ functions with support $U \cup S$
 and $V \cup  S$ respectively, and let $\Omega= (f- g)dx_1\wedge 
\ldots\wedge  dx_n$. 
Then $S$ is isoperimetric for prescribed $\Omega$-volume; indeed it is the 
only surface with its $\Omega$-volume.

A similar hypothesis appears for example in [\cite{DF2}, Thm.\ 5.1].
\\[4mm]
\noindent
\em Proof of Proposition $2.5.$ \em 
 For every smooth variation vectorfield $\mathbf{v}$  on $S$, initially
$$\mathrm{d}V\!\!/\!\mathrm{d}\,t = \int_S 
(\mathrm{d}\omega \LA{\llcorner} T)(\mathbf{v}),$$
basically because by Stokes's theorem the change in volume is the integral 
of $\mathrm{d} \omega$ over the volume swept out.
It follows immediately that if the mean curvature vector $\mathbf{H}$
 is proportional to $\mathrm{d}\omega \LA{\llcorner} T$ (weakly), i.e. if
$$\mathrm{d}A\!/\!\mathrm{d}\,t = \lambda \int_S 
(\mathrm{d}\omega \LA{\llcorner} T)(\mathbf{v}),$$
then $V$ constant implies that $\mathrm{d}A\!/\!\mathrm{d}\,t = 0$, 
so $S$ is stationary.

Conversely, suppose that $S$ is stationary and 
$\mathrm{d}\omega  \LA{\llcorner} T$
 is not identically 0. Then the constraint is nonsingular as well as smooth, 
so for some Lagrange multiplier  ${\lambda}$, 
$\mathrm{d}A\!/\!\mathrm{d}\,t 
= \lambda (\mathrm{d}V\!\!/\!\mathrm{d}\,t)$, 
as desired. Alternatively suppose that $S$ is stationary, that $M = 
\mathbf{R}^n$, that $\mathrm{d}\omega$ is constant, and that 
$\mathrm{d}\omega  \LA{\llcorner} T
 $ is 0 almost everywhere. Then variations of the form $S_t = S + t\mathbf{v}$
with $\mathbf{v}$ of the special form $\mathbf{v} = \varphi\cdot \mathbf{v}_0$
for some smooth scalar function $\varphi$ and fixed vector $\mathbf{v}_0$
preserve $\omega$-volume. Since $S$ is stationary, 
$\mathrm{d}A\!/\!\mathrm{d}\,t$ is initially
0. Since such vectorfields $\mathbf{v}$ span the space of all smooth variation 
vectorfields, $\mathbf{H}$ is 0.

Finally suppose that $M = \mathbf{R}^n$, $\mathrm{d}\omega$ is constant and 
simple, 
and $S$ is stationary, smooth, connected, and bounded. We may assume that 
$\mathrm{d}\omega$  is $\mathrm{d}x_1\wedge \ldots \wedge \mathrm{d}x_{m+1}$. 
For coordinates $(x, y)$
 on $\mathbf{R}^{m+1}\times \mathbf{R}^{n-m-1}$,  consider a family of 
diffeomorphisms given on $S$ by $F_t(x, y) = (x, y/(1+t))$. Since they 
preserve volume, $\mathrm{d}A\!/\!\mathrm{d}\,t$ initially must be 0, 
which means 
that $S$ is 
everywhere horizontal and lies in a horizontal copy of $\mathbf{R}^{m+1}$. 
Since $S$ has constant mean curvature (nonzero because $S$ is bounded), 
$S$ is round by Alexandrov's Theorem.\\[5mm]
\noindent
Finally we consider the third case of prescribed multi-volume.
\begin{proposition}[\cite{M4}, Thm.\ 2.2]
A boundary $S$ with unit tangent 
$m$-vector  $T$ in $\mathbf{R}^n$ is stationary for prescribed multi-volume 
if and only if for some constant $(m+1)$-form $\Omega$, the mean curvature 
$\mathbf{H}$ of $S$ weakly satisfies 
$m\,\mathbf{H}\cdot = \Omega \LA{\llcorner} T$, 
i.e., for any smooth variation vectorfield $\mathbf{v}$, initially
\begin{equation}
 \mathrm{d}A\!/\!\mathrm{d}\,t = -\int_S m\, \mathbf{H}\cdot \mathbf{v} = 
-\int_S (\Omega \LA{\llcorner} T)(\mathbf{v}) 
= -\int_S \Omega(T\!\wedge \mathbf{v}) .
\end{equation}
\end{proposition}
\noindent
\em Proof. \em Assume that (2) holds. Choose a smooth $\omega$ such that 
$\Omega = \mathrm{d} \omega$. By Proposition 2.5, $S$ is stationary for 
prescribed 
$\omega$-volume, and hence for prescribed multi-volume. 

Conversely, assume 
that $S$ is stationary. For any covector $\Omega = \mathrm{d} \omega$, consider 
the associated volume $\int_S \omega$. (Fixing multi-volume is equivalent 
to fixing all such volumes or just the axis volumes $V_I$.) 
For every smooth variation vectorfield $\mathbf{v}$ on $S$, initially
$$\mathrm{d}V_I\!/\!\mathrm{d}\,t = \int_S(dx_I \LA{\llcorner}  T)
(\mathbf{v}).$$
We consider variations of the forms $S_t= S+ t\mathbf{v}$ with $\mathbf{v}$ 
of the special form $\mathbf{v} = \varphi\cdot \mathbf{v}_0$ for some smooth 
scalar function $\varphi$ and fixed vector $\mathbf{v}_0$, which span the space 
of all smooth variations. A variation of this simple form never alters 
volumes for $\Omega$ outside
$$
\mathrm{span}\{T\!\wedge \mathbf{v}_0 :
\mathbf{v}_0\in \mathbf{R}^n, \mbox{ values of } T \mbox{ at Lebesgue points}
\} \subset  \wedge^{m+1}\mathbf{R}^n.
$$
The constraint for such volumes is nonsingular as well as smooth, as can be 
seen by consideration of variations supported in small neighborhoods of 
Lebesgue points of $T$. Therefore for some Lagrange multiplier 
 $\boldsymbol{\lambda}=(\lambda_I)$, 
$\mathrm{d}A\!/\!\mathrm{d}\,t = \sum \lambda_I  
(\mathrm{d}V_I\! /\!\mathrm{d}\,t)$, 
so with 
$\Omega = -\sum \lambda_I dx_I$
$$ \mathrm{d}A\!/\!\mathrm{d}\,t 
= - \int_S (\Omega \LA{\llcorner} T)(\mathbf{v}),$$
as desired.

\section{Existence and Regularity of Isoperimetric Surfaces}
\noindent
This section presents standard geometric measure theory results on existence and
regularity.\\[2mm]
{\bf 3.1.~Existence.}
\em 
 If $M$ is compact or if $M = \mathbf{R}^n$ and $\mathrm{d}\omega$ 
is constant, then 
isoperimetric surfaces $S$ exist for all prescribed volumes $\mathrm{v}(S)$, 
$\omega$-volumes, and multi-volumes and are compact.\em
\\[3mm]
\noindent
\em Proof. \em  If $M$ is compact  there are no issues, one just takes a 
minimizing sequence and applies the Compactness Theorem [\cite{M1}, 
Chap.\ 5] to get 
a solution in the limit. In $\mathbf{R}^n$, local compactness still provides 
a possibly unbounded area-minimizing limit among locally integral currents 
[\cite{M1}, \S 9.1]. By Propositions 2.3, 2.5, 2.6, an 
isoperimetric surface has 
constant-magnitude mean curvature. Thence ``monotonicity'' [\cite{A1}, 5.1(3)]
yields a positive lower bound on the area inside a unit ball about every 
point of $S$, and it follows that $S$ is compact. 

A more serious problem is  that there may be volume loss to infinity. 
One uses a concentration lemma and  translation to obtain a minimizer with 
nonzero volume [\cite{M1}, \S 13.4]. Then for  prescribed volume 
$\mathrm{v}(S)$ or 
prescribed $\omega$-volume, one uses scaling  (and a flip of orientation 
if necessary) to obtain the prescribed volume.  For prescribed multi-volume 
one repeats the process countably many times to  recover all the volume 
[\cite{M1}, \S 13.4].\\[4mm]
\noindent
{\bf 3.2.~Regularity.} By Allard's regularity theorem [\cite{A1}, Sect.\ 8], 
any surface with weakly bounded mean curvature is a $C^{1,\alpha}$ 
submanifold on 
an open dense set. By Propositions 2.3, 2.5, 2.6, this includes all three 
types of isoperimetric surfaces in $\mathbf{R}^n$, assuming $\mathrm{d}\omega$
constant. It probably includes isoperimetric surfaces for prescribed volume 
$\mathrm{v}(S)$ in smooth Riemannian manifolds, but we do not know how to prove that.

For a negative example for prescribed $\omega$-volume in $\mathbf{R}^n$ with 
$\mathrm{d}\omega$  nonconstant, see  the Remark after Proposition 2.5.

It is not known whether isoperimetric surfaces for prescribed volume 
$\mathrm{v}(S)$ in
smooth Riemannian manifolds and for prescribed multivolume in $\mathbf{R}^n$
 enjoy the same regularity as area-minimizing surfaces without volume 
constraints, even for the easier Lagrange multiplier problem; cf.\ 
[\cite{M1},  Chap.\ 8], [\cite{DF1}, \S 5], [\cite{DS}, Intro.\  
and 5.5(iii)].  In general, it is not 
known  even whether a tangent cone is minimizing, because the cost of small 
volume  adjustments is not known to be linear. (Note e.g.\ the extra hypothesis
required in [\cite{DF2}, Thm.\ 5.1].)

\section{Round Spheres Uniquely Minimizing or Stable?}
\noindent
This section proves that round spheres are uniquely minimizing for all three 
volume constraints and conjectures that they are uniquely stable in 
$\mathbf{R}^n$ for prescribed volume $\mathrm{v}(S)$ and for prescribed 
$\Omega$-volume  for $\Omega$ constant (but not for prescribed multi-volume).

\begin{proposition} Round $m$-spheres $S_0$ are uniquely minimizing for all 
three cases; for prescribed $\omega$-volume (case 2) we need to assume 
$\Omega =\mathrm{d}\omega$ constant and maximum on the $(m+1)$-ball 
bounded by $S_0$.
\end{proposition}
\noindent
\em Proof.  Case $(1)$, prescribed volume $\mathrm{v}(S)$. \em 
Almgren \cite{Alm1}, indeed $\mathrm{mod}\,\nu$  for all $\nu$.\\[2mm]
\em Case $(2)$, prescribed $\omega$-volume. \em 
 We may assume $\Omega =\mathrm{d}\omega$  is 1 on the disc $D$. Now let $S$ be
any surface with the same $\omega$-volume, and let $R$ be a volume-minimizing 
surface bounded by $S$. Then
$$ |D| = \int_D \Omega  = \int_R \Omega \leq  |R|.$$
By Case (1), $|S_0| \leq  |S|$, with equality only if $S$ is a round sphere 
and $\Omega = 1$ on $R$.\\[2mm]
\em Case $(3)$, prescribed multi-volume, \em  follows from Case (2) with 
$\Omega$  the simple form dual to
the disc.\\[4mm]
\noindent
\em Remark. \em Salavessa \cite{S2} proves the weaker result that 
associated round 
spheres have nonnegative second variation for prescribed $\Omega$-volume 
for $\Omega$ the K\"{a}hler form on $\mathbf{R}^6$.\\[4mm]
\noindent
 The following conjecture 
would generalize a codimension-1 stability theorem of Barbosa
and do Carmo \cite{BdC} to higher codimension.\\[4mm]
{\bf  4.2. Conjecture.} \em  In $\mathbf{R}^n$, round $m$-spheres $S_0$
 are the only smooth stable surfaces $S$ for given volume $\mathrm{v}(S)$
 or given $\Omega$-volume for constant $\Omega$ $($although not for given
 multi-volume \em [\cite{M4}, Cor.\ 3.2]). \\[3mm]
\noindent
\em Proof for $\Omega$-volume for $m=1$. \em 
By [\cite{M4}, Thm.\ 3.1], which applies to stationary as well as 
minimizing curves, 
a stationary closed curve for prescribed multi-volume or equivalently
for prescribed $\Omega$-volume is of the form
$$ C(s) = a_0 + a_1e^{iw_1s}e_1 +\ldots + a_k e^{iw_ks}e_{2k-1},$$
with $w_j$ increasing positive integers. If $k = 1$, this curve is a circle. 
If $k > 1$, this curve is neither minimizing nor stable for given 
$\Omega$-volume: in the second component, which encircles the origin twice, 
enlarging one loop and shrinking the other reduces length to second order 
for fixed area. (This variation does not preserve multi-volume because it
alters area in the $e_{14}$ plane for example. Indeed, this curve is 
minimizing for prescribed multi-volume [\cite{M4}, Cor.\ 3.2].)\\[4mm]
\em Remarks. \em 
 The proof that round spheres are the only minimizers dates from 1986
[Almgren \cite{Alm1}]. Generalizing the codimension 1 proof of 
Barbosa-do Carmo \cite{BdC}
and Wente \cite{W}
 seems to need mean curvature parallel on $S$ (see \S 4.3 below). 
Otherwise second variation in the normal direction is more positive 
([\cite{Sc}, p.\ 171] or because $\mathrm{d}H\!/\mathrm{d}\,t$  
is less negative). 
It also seems to need 
normal bundle geometrically trivial. Salavessa \cite{S1} seems 
further to need a 
Minkowski-type hypothesis: in her completely different
terminology ``$\int_M S(2 + h\|H\|)dM \leq  0.$''
If the scalar mean curvature $H$ is not constant, one
could replace $H$ with its average except that the average value of the 
square of $H$ is greater than the square of the average value. 
For prescribed volume $\mathrm{v}(S)$, this conjecture remains open even for 
curves $(m=1)$ in $\mathbf{R}^3$; perhaps the circle is even the only
isoperimetric-stationary one-component curve in $\mathbf{R}^n$.

Suppose that $C$ is a counterexample for $m=1$ in $\mathbf{R}^3$ for say 
volume $\pi$. Further suppose that $C$ has as expected (see 2.2 and 2.3) 
curvature $\kappa$ of magnitude $|C|/2\pi > 1$ in
the direction of the inward normal. 
By the isoperimetric inequality (or by Bol-Fiala for
a disc), $|C| > 2\pi$. Alternatively, since $\kappa = |C|/2\pi$, 
$|C| = 2\pi \kappa$, Gauss-Bonnet yields $C^2/2\pi= \int |\kappa |> 2\pi\chi$ , 
again yielding $|C| > 2\pi$  for a disc, 
since the Euler characteristic $\chi \leq  1$ for
cases of interest ($R$ connected). Moreover, since the curvature of $R$
 along $C$ vanishes, so
does the curvature of $R$ normal to $C$: locally as a graph 
$f_{xx}$ and $f_{yy}$ vanish, but not necessarily $f_{xy}$, 
so we don't see how to prove e.g.\ that $f_{yyy}$ vanishes and that $R$
 contains rays from the boundary and must be a flat disk.\\[4mm]

\noindent
{\bf 4.3. Second Variation. } The formula for the second variation, 
that is, the second derivative of area for a smooth family of perturbations, 
is given by Schoen [\cite{Sc}, p.\ 171]. Note that every variation 
vectorfield for a 
compact surface $S$ for prescribed constant $\Omega$-volume in $\mathbf{R}^n$
with $\mathrm{d}V\!\!/\!\mathrm{d}\,t$ initially 0 is part of a 1-parameter 
family with fixed volume 
obtained by adjusting any 1-parameter family with 
$\mathrm{d}V\!\!/\!\mathrm{d}\,t$ initially 0 by 
continuous rescalings by homotheties. This generalizes to exact nonconstant 
$\Omega$-volume in manifolds as long as there is a variation vectorfield 
(like the one for scaling) for which $\mathrm{d}V\!\!/\!\mathrm{d}\,t$ 
is not zero. This
corresponds to the fact that if $f$ is a smooth function on $\mathbf{R}^n$, 
$ \partial\! f\!/\partial x_1 = 0$, and $\mathrm{grad} f \neq  0$, then $f$
vanishes on a smooth horizontal curve through 0.

On the other hand, even minimizers for prescribed multi-volume can be unstable
for variations which preserve multi-volume to first order (and hence cannot 
correspond to 1-parameter multi-volume-preserving families). For example, 
consider the curve $(e^{is}, e^{2is})$
in $\mathbf{R}^2\times \mathbf{R}^2$, which is isoperimetric [\cite{M4}, 
Cor.\ 3.2]. 
In the second factor the curve is two copies of the unit circle. 
Shrinking one and expanding the other preserves multi-volume
to first order but reduces length to second order. 

Conjecture 2.2(2) 
would imply that this curve is not even stationary for prescribed
volume $\mathrm{v}(S)$. The unique volume-minimizing surface bounded by this 
curve is  $\{w = z^2\}$, because complex analytic varieties are uniquely volume 
minimizing [\cite{M1}, 6.3]. Note that $\mathbf{H}$ is proportional to 
$-(e^{is}, 4e^{2is})$, while the inward conormal is proportional to 
$-(\mathrm{d}z, 2z\mathrm{d}z)$ hence to $-(e^{is}, 2e^{2is})$. 
By Proposition 2.2 and Remarks, the curve  is not  stationary.

Incidentally, this curve is not a graph over every axis plane: 
its projection in the $x_1x_4$-plane is  a figure 8 enclosing  signed area $0$.
In this case the problem of prescribing \em unsigned \em 
areas would have a different 
solution, presumably circles in the $x_1x_2$-plane and the 
$x_3x_4$-plane.\\[2mm]

\noindent
{\bf 4.4. Yau \cite{Y} on parallel mean curvature vector}. 
Yau \cite{Y}  proved that 
every smooth 2-dimensional 
 surface in $\mathbf{R}^n$ with parallel mean curvature vector 
is one of four types:\\[3mm]
(1) constant-mean-curvature hypersurfaces in some $\mathbf{R}^3 \subset  
\mathbf{R}^n$,\\[2mm]
(2) constant-mean-curvature hypersurfaces of some 
$\mathbf{S}^3 \subset \mathbf{R}^n$,\\[2mm]
(3) minimal submanifolds of some hypersphere $\mathbf{S}^{n- 1} \subset 
\mathbf{R}^n$,\\[2mm]
(4) minimal submanifolds of $\mathbf{R}^n$.\\[2mm]

One could study Lawson \cite{Law} and other examples of minimal surfaces in 
$\mathbf{S}^3 \subset  \mathbf{R}^4$ (type (2) and (3)).\\[4mm]
\noindent
\em Example. \em $S = \mathbf{S}^1\times  \mathbf{S}^1$ in $\mathbf{R}^4$. 
$\mathbf{H}$ is  parallel  (Yau type (2) and (3)), but $S$ is not stationary
even for fixed multi-volume, because it and all of its scalings have 
multi-volume 0. In particular there is no constant 3-form 
$\Omega$  dual to $T\!\wedge \mathbf{H}$  on $S$, a completely trivial
consequence, since codimension-1 forms are simple and $T\!\wedge \mathbf{H}$
 is not constant. An obvious variable calibration candidate, 
$\mathrm{d}r \,\,  r_1\mathrm{d}\theta _1 \, r_2\mathrm{d} \theta_2$ 
is not closed. We think that for general
reasons there is a smooth classical calibration $\Omega$ of a small band of 
the cone over $\mathbf{S}^1\times \mathbf{S}^1$,
as for any small stationary surface (Lawlor \cite{L}), and $S$ is stationary 
for prescribed $\Omega$-volume. It is the same story for any stationary 
product of spheres or of minimal submanifolds of spheres.

\section{Calibrations }

The classical theory of calibrations (see [\cite{M1}, \S 6.4] and 
references therein) 
says that if there is a closed form $\omega$ on a smooth Riemannian manifold 
such that $|\omega| \leq  1$ with equality on the tangent planes to a surface 
$S$, then $S$ is area minimizing in its homology class. The form $\omega$
 is called a calibration of $S$. Morgan \cite{M5} noted that for hypersurfaces
if the condition that $\mathrm{d}\omega$ be 0 is relaxed to the condition that 
$\mathrm{d}\omega$ be a constant multiple of the volume form, then 
$S$ still minimizes 
area for prescribed volume. In particular, constant-mean-curvature graphs have 
such ``$\mathrm{d}$-constant calibrations'' \cite{M5}, citing 
[\cite{M1}, \S 6.1].

The following proposition is a trivial extension to prescribed 
$\omega$-volume in general codimension.

\begin{proposition} If $\omega$ attains its maximum value (say 1) 
everywhere on a surface $S$ in a smooth Riemannian manifold, then $S$
 is isoperimetric for prescribed $\omega$-volume.
\end{proposition}

\noindent
\em Remarks. \em  In particular, every smooth surface is isoperimetric for 
some smooth $\omega$.

By Proposition 2.5, unless $\mathrm{d}\omega \LA{\llcorner}  T$  
vanishes on every unit 
tangent plane $T$ to $S$, the mean curvature vector $\mathbf{H}$ of $S$
 is proportional to $\mathrm{d}\omega \LA{\llcorner}  T$.\\[4mm]
\noindent
\em Proof of Proposition $5.1$. \em If $S'$ has the same $\omega$-volume, then
$$ |S| = \int_S \omega = \int_{S'} \omega \leq  |S' |. $$\\[-3mm]

\noindent
\em Remarks. \em Consider a smooth exact $(m+1)$-form $\Omega$ and a surface 
$S$ with unit tangent planes $T$ and mean curvature vector $\mathbf{H}$ such 
that $\mathbf{H} = \Omega \LA{\llcorner} T$. By Proposition 2.5, $S$ is 
stationary for given $\Omega$-volume. By Proposition 5.1,  
 $S$ is area minimizing for 
given $\Omega$-volume. This probably always holds locally  \em a  la \em 
Lawlor \cite{L}. Conversely, if $S$ is area minimizing for given 
$\Omega$-volume,  there is probably in some generalized weak sense a 
calibration $\omega$ \em a  la \em  Federer \cite{F}.

\end{document}